\documentclass{amsproc}
\usepackage{amssymb}
\usepackage{graphicx}
\usepackage{amscd}
\usepackage{amsmath}
\newcommand{\CI}{\mathcal{I}}
\newcommand{\CJ}{\mathcal{J}}

\theoremstyle{plain}

\newtheorem{definition}{Definition}

\newtheorem{remark}{Remark}

\newtheorem{theorem}{Theorem}
\numberwithin{equation}{section}

\newcommand\bess{\begin{eqnarray*}}
\newcommand\eess{\end{eqnarray*}}

\newcommand\beq{\begin{equation}}
\newcommand\eeq{\end{equation}}

\def\Nset{\mbox{I\kern-.21em N}}
\def\RE{{\mbox{\rm I\kern-.21em R}}}
\def\ZZ{{\mbox{\sf Z\kern-.45em Z}}}
\def\vv{\kern.344em{\rule[.18ex]{.075em}{1.32ex}}\kern-.344em}

 \def\gm{\gamma}  
  \def\Om{\Omega}

\def\<{\langle} \def\>{\rangle}

\def\tr{\operatorname{tr}}

\begin{document}

\title[]{On an inverse problem for tree-like networks of elastic strings}
\date{September, 2008}
\author{Sergei Avdonin$^1$, G\"unter Leugering$^2$ and Victor Mikhaylov$^1$}
\address{$^1$Department of Mathematics and Statistics \\
University of Alaska Fairbanks\\
PO Box 756660\\
Fairbanks, AK 99775} \address{$^2$Department Mathematik\\
Lehrstuhl Angewandte Mathematik II\\
Friedrich-Alexander-University Erlangen-N\"urnberg\\
Martenstr. 3\\
D 91058 Erlangen} \email{ffsaa@uaf.edu,
leugering@am.uni-erlangen.de, ftvsm@uaf.edu} 
 \maketitle

\begin{abstract}
We consider the in-plane motion of elastic strings on tree-like
network, observed from the 'leaves'. We investigate the inverse
problem of recovering not only the physical properties i.e. the
'optical lengths' of each string, but also the topology of the tree
which is represented by the edge degrees and the angles between
branching edges. To this end use the boundary control method for
wave equations established in~\cite{AK,B}. It is shown that under
generic assumptions the inverse problem can be solved by applying
measurements at all leaves, the root of the tree being fixed.
\end{abstract}

\section{Introduction}

\section{Forward dynamical and spectral problems for the two-velocity system on the tree.}
In many problems in science and engineering network-like
structures play a fundamental role. The most classical area of
applications consists of flexible structures made of strings,
beams, cable and struts. Bridges, space-structures, antennas,
transmission-line posts, steel-grid structures as reinforcements
of buildings and other projects in civil engineering. See Lagnese,
Leugering and Schmidt~\cite{LLS1994} for an account of multi-link
structures. More recently applications also on a much smaller
scale came into focus. In particular hierarchical materials like
ceramic or metallic foams, percolation networks and even carbon
nano-tubes have attracted much attention. In the latter context,
the problem is understood as a quantum-tree-problem. See e.g.
Kuchment~\cite{kuch}, Kostrykin and Schrader~\cite{ks1}, Avdonin
and Kurasov~\cite{AK}. In all of these areas the topology of the
underlying networks or graphs plays a dominant role. The
understanding of the influence of the local topology and physical
parameters, say at a given branching point, on the global
mechanical or scattering properties is crucial in this area.
Failure detection in mechanical multi-link structures by
non-invasive methods as well as topological and material
sensitivities with respect to an observer play an important role.
Gaining this understanding is the central focus of this paper.
Undoubtedly, the inverse problem for mechanical structural
elements like membranes and plates has been discussed in the
literature. The famous question by Kac~\cite{Kac}: "can one hear
the shape of a drum" initiated major research in this direction.
This question has been repeated in the literature regarding other
structures, also for string-networks on a tree by Avdonin and
Kurasov~\cite{AK}. However, in that work strings have been
considered as deflecting out of the plane rather than in the
plane. The important and in fact crucial difference is that such
networks are insensitive for the topology of the graph in the
sense that the coupling conditions do not reflect the angles at
which the strings are 'glued' together. Only in the case of
in-plane motion are the coupling conditions dependent on the local
geometry of the multiple joints. This observation is even more
relevant for networks containing beams, a case that is subject to
current investigation. In a more abstract setting, where also
electromagnetic or quantum effects are considered on graphs, one
observes that such an in-plane modeling involves multi-channel and
multi-velocity models for wave propagation in thin structures.

Tackling inverse problems involves the understanding
Steklov-Poincar{\'e}-operators just as in the case of domain
decomposition. Such operators for problems on graphs have been
investigated in Lagnese and Leugering~\cite{LagneseLeugering2004}.
Scattering matrices, indeed the Tichmarsh-Weyl function for
in-plane-networks of strings, at that time called echo-analysis,
have been investigated in Leugering~\cite{Leugering1996}.  In
particular, the understanding of controllability properties of the
underlying structures is crucial for a dynamic, and indeed
real-time, detection of physical and geometrical properties.
Again, exact controllability of networks of strings both in the
out-of-the-plane and the more important in-plane mode has been
investigated by Lagnese, Leugering and
Schmidt~\cite{LagneseLeugeringSchmidt1993,LagneseLeugeringSchmidt1994,LLS1994},
see also Avdonin and Ivanov~\cite{AvdoninIvanov1995}. There it has
been shown that under generic assumptions, controllability of a
rooted tree holds by controls at the leaves. Later Zuazua and
Leugering~\cite{LeugeringZuazua1999} showed that under more
refined assumptions on the nature of the out-of-the-plain
string-tree exact controllability in refined spaces was even
possible when the root was controlled only. This research has been
extended considerably in Dager and Zuazua~\cite{DagerZuazua2006}.
As it turned out in Belishev~\cite{B,BV} Avdonin and
Kurasov~\cite{AK} exact controllability of both the state and the
velocity appeared too demanding. Indeed, their
'boundary-control-approach' is based on controllability of the
state only. The work on inverse problems by the way of the
boundary-control-approach has by now become a major tool.  The
current paper is no exception in that direction. For inverse
problems on graphs see also the work of Yurko \cite{Y}.

\begin{figure}
\centering
\includegraphics[width=2in]{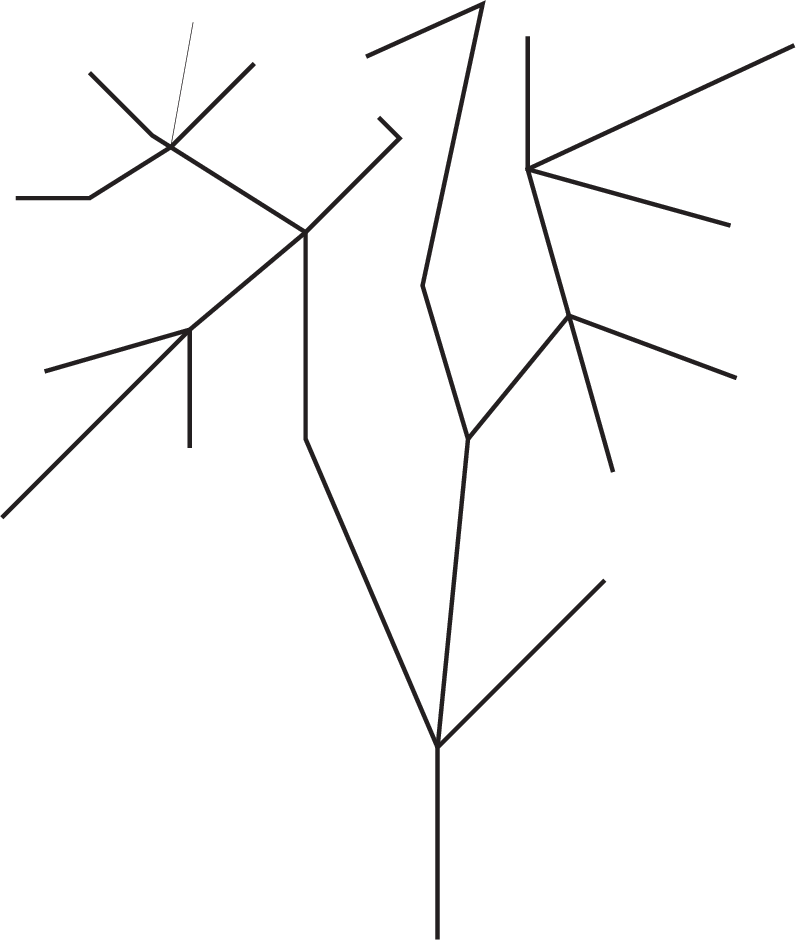}
\caption{The tree}\label{fig2}
\end{figure}

Let $\Omega$ be a finite connected compact planar graph without
cycles, i.e. a tree. The graph consists of edges
$E=\{e_1,\ldots,e_N\}$ connected at the vertices
$V=\{v_1\ldots,v_{N+1}\}$, see figure \ref{fig2}.

Every edge $e_j\in E$ is identified with an interval
$(a_{2j-1},a_{2j})$ of the real line. The edges are connected at the
vertices $v_j$ which can be considered as equivalence classes of the
edge end points $\{a_j\}$. The boundary
$\Gamma=\{\gamma_1,\ldots,\gamma_m\}$ of $\Omega$ is a set of
vertices having multiplicity one (the exterior nodes). We suppose
that the each edge of the graph is equipped with the densities
\begin{equation}
\label{density} \frac{1}{k^2_{i1}},\,\frac{1}{k^2_{i2}} \quad
\text{on $e_i$},\ i=1\ldots,N.
\end{equation}
The densities on the edges determine the density on the graph by
the rule: for $x\in e_i\subset \Omega\backslash V$, $\rho(x)$ is a
constant which is equal to
$\max{\left\{\frac{1}{k^2_{i1}},\,\frac{1}{k^2_{i2}}\right\}}$.

\begin{figure}
\centering
\includegraphics[width=1.5in]{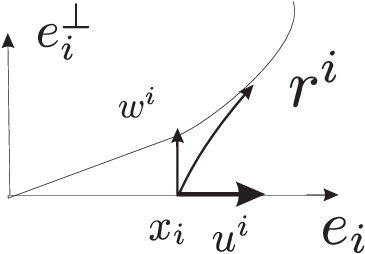}
\caption{Representation of planar displacement}\label{fig1}
\end{figure}

Since the graph under consideration is a tree, for every
$a,b\in\Omega,$ $a\not=b,$ there exist the unique path $\pi[a,b]$
connecting these points. The density determines the optical metric
and the optical distance
\begin{eqnarray*}
d\,\sigma^2:=\rho(x)|d\,x|^2,\quad x\in\Omega\backslash V,\\
\sigma(a,b)=\int\limits_{\pi[a,\,b]}\sqrt{\rho(x)}|d\,x|, \quad
a,b\in\Omega,
\end{eqnarray*}
The optical diameter of the graph $\Omega$ is defined as
\begin{equation*}
\label{diam} d(\Omega)=\max_{a,\,b\in\Gamma}\sigma(a,b).
\end{equation*}
The graph $\Omega$ and the optical metric determine the {\it
metric graph} denoted by $\{\Omega,\rho\}$. For a rigorous
definition of the metric graph see, e.g.
\cite{ks1,kuch,KuSt,AM1,graphbook}. 
The densities can also be considered as reciprocals of a stiffness
parameter $k_{i1}^2,k_{i2}^2$ related to a horizontal and vertical
displacement $u^i,\, w^i$ of an elastic string along the edge with
index $i$. Each edge carries along an individual coordinate system
$e_i,\; e_i^\perp$ such that the displacement
\[ r^i(x):= u^i(x) e_i + w^i(x) e_i^\perp\]
splits into a longitudinal displacement $u^i(x)$ and a vertical
displacement $w^i(x)$ at the material point $x\in
[a_{2i-1},a_{2i}]$. See figure \ref{fig1}.

The stiffness matrix can then be expressed as
\[
K_i= k_{i1}^2 e_i e_i^T+k_{i2}^2 e_i^\perp (e_i^\perp)^T.\] For the
sake of completeness we describe in more detail the full set of
equations in more generality than we actually need later. Given a
node $v_J$ we define
\[\CI_J := \{i\in\CI | e_i \text{ is incident at } v_J\}\]
the incidence set, and $d_J=|\CI_J|$ the edge degree of $v_J$. The
set $\CJ=\{J| v_J\in V\}$ of node indices splits into simple nodes
$\CJ_S$ and multiple nodes $\CJ_M$ according to $d_J=1$ and $d_J>1$,
respectively. We introduce the subsets of boundary nodes according
to tangential and vertical motions:
\[ \CJ_D^t:=\{ J\in\CJ_S | r^i(v_J)\cdot e_i=0\}\]
\[ \CJ_D^n:=\{ J\in\CJ_S | r^i(v_J)\cdot e_i^\perp=0\}\]
\[ \CJ_N^t:=\{ J\in \CJ_S | d_{iJ}K_ir^i_x(v_J)\cdot e_i=0\}\]
\[ \CJ_N^n:=\{ J\in \CJ_S | d_{iJ}K_ir^i_x(v_J)\cdot e_i^\perp=0\},\]
where
\[ d_{iJ}=\begin{cases} -1 & \text{ edge i starts at } v_J\\
1 & \text{ edge i ends at } v_J \end{cases},\;\;
x_{iJ}=\begin{cases} 0& d_{iJ}=-1\\\ell_i & d_{iJ}=1.\end{cases}\]
and $r^i(v_J):=r^i(x_{iJ}$.

Notice that these sets are not necessarily disjoint. Obviously, the
set of completely clamped vertices can be expressed as
\begin{equation}\label{clamped}
\CJ_D^0:= \CJ_D^t\cap\CJ_D^n
\end{equation}
Similarly, a vertex with completely homogenous Neumann conditions is
expressed as $\CJ_N^n\cap\CJ_N^t$. At tangential Dirchlet nodes in
$\CJ_D^t$ we may, however, consider normal Neumann-conditions as in
$\CJ_N^n$ and so on. We may then consider time dependent
displacements $r^i(x,t)$, $t\in[0,T]$. The system of equations
governing the full transient motion is then given by
\begin{equation}\label{dynamic}
\left\{ \begin{gathered}
r^i_{tt}-K_i r^i_{xx}+ c_i r^i = f^i\quad \in (0,\ell_i)\\
r^i(v_D)= u_D, \; i\in \mathcal{I}_D, D\in\mathcal{J}_D\\
  d_{iJ}r^i_x(v_N)=g_J\quad i\in \mathcal{I}_N, N\in\mathcal{J}_N\; \\
r^i(v_J)=r^j(v_J)\quad i,j\in \mathcal{I}_J, J\in\mathcal{J}_M\; \\
\sum\limits_{i\in\mathcal{I}_J}d_{iJ}K_i r^i_x(v_J)=0\quad
J\in \mathcal{J}_M,\\
r^i(\cdot,0)=r^i_0,\; {r}^i_t(\cdot, 0)=r^i_1
\end{gathered}\right.
\end{equation}
In this representation we used capital letters for vertices in order
to improve the readability of the formulae. It important to
understand the coupling conditions \eqref{dynamic}$_{4,5}$. Indeed,
the first of these conditions simply expresses the continuity of
displacements across the vertex $v_J$. Without this condition the
network falls apart. The second condition reflects the physical law
that the forces at the vertex $v_J$, in the absence of additional
external forces acting on $v_J$, should add up to zero. Notice that
the coupling at multiple nodes $v_J$, those where $|I_J|>1$, is a
vectorial equation. This is in contrast to the out-of-the-plane
model, where no such vectorial couplings occur which, in turn, makes
the problem then independent of the particular geometry. In the case
treated here the geometry, represented by the pairs $(e_i,
e_i^\perp)$ does play a crucial role.

In Leugering and Sokolowski~\cite{LeugeringSokolowski2009} the
static system i.e. \eqref{dynamic} without time-dependence,  has
been investigated with respect to topological sensitivities, such
that the potential energy or other functionals are considered under
variation of the edge-degree, nodal positions and edge-deletion.
Like in this paper the analysis of the Steklov-Poincar{\'e}
operators play a crucial role.

Another  remark about model \eqref{dynamic} is in order. If the
system is static, the stiffness only active for longitudinal
displacements, and if the state is edge-wise linear, then
\eqref{dynamic} comes down to a truss-model.

In the global Cartesian coordinate system one can represent each
edge by a rotation matrix
\[
R_{\alpha_i}=\begin{pmatrix} \cos{\alpha_i} &
-\sin{\alpha_i}\\
\sin{\alpha_i} & \cos{\alpha_i}
\end{pmatrix}
\]
where $e_i=(\cos \alpha_i, \sin \alpha_i)^T$. In fact, it is the
global coordinate system that we will use throughout the paper, as
we are going to identify the angles $\alpha_{ij}$ between two
branching edges.

\subsection{Spectral settings}

Let us introduce the spaces of functions on the graph $\Omega$:
\begin{equation}
L_2(\Omega)=\bigoplus_{i=1}^N L_2(e_i, \mathbb{R}^2)
\end{equation}
For the element $U\in L_2(\Omega)$ we write
\begin{equation}
U=\left\{u, w\right\}=\left\{
\begin{pmatrix}u^i\\
w^i
\end{pmatrix}
\right\}_{i=1}^N,\quad u^i,w^i\in L_2(e_i).
\end{equation}

Let us consider the internal vertex $v\in V\backslash\Gamma$ and
the edges $e_1,\ldots e_K$ incident at $v$. We denote by
$\alpha_{ij}$ the angle between two edges $e_i$ and $e_j$ counting
from $e_i$ counterclockwise. Let $R_\alpha$ be the matrix
corresponding to the rotation on the angle~$\alpha$:
\begin{equation}
R_{\alpha}=\begin{pmatrix} \cos{\alpha} &
-\sin{\alpha}\\
\sin{\alpha} & \cos{\alpha}
\end{pmatrix}
\end{equation}
Let us introduce the matrices
\begin{equation}
D'_i=\begin{pmatrix}
k_{i1}^2 & 0\\
0 & k_{i2}^2
\end{pmatrix},\quad i=1,\ldots N.
\end{equation}

We can now reformulate the compatibility conditions in
\eqref{dynamic} at multiple nodes (vertices) $v$ using global
coordinates. For the sake of self-consistency in this framework we
put this in the format of definitions.

\begin{definition}
We say that the edge-wise continuous function $U$ satisfy the first
condition (continuity) at the (multiple node) internal vertex $v$ if
\begin{equation}
\label{C_ContInV}
\begin{pmatrix}
u^i(v)\\
w^i(v)
\end{pmatrix}=R_{\alpha_{ij}}\begin{pmatrix}
u^j(v)\\
w^j(v)
\end{pmatrix},\quad 1\leqslant i,j\leqslant K,\,\,\, i\not=j.
\end{equation}
\end{definition}
\begin{definition}
We say that the edge-wise continuously differentiable function $U$
satisfies the second condition (force balance) at the internal
vertex $v$ if
\begin{equation}
\label{C_KirhInV} D'_1\begin{pmatrix}
u^1_x(v)\\
w^1_x(v)
\end{pmatrix}=\sum_{j=2}^K R_{\alpha_{1j}}D'_j\begin{pmatrix}
u^j_x(v)\\
w^j_x(v)
\end{pmatrix}.
\end{equation}
\end{definition}

We associate the following spectral problem to the graph:
\begin{eqnarray}
\left.\begin{array}r -k_{i1}^2\phi_{xx}^i=\lambda \phi^i\\
-k_{i2}^2\psi_{xx}^i=\lambda \psi^i
\end{array}\right\}\quad x\in e_i,\quad i=1\ldots,N, \label{SpEq1ch}\\
\text{$\left\{\phi,\psi\right\}$ satisfies (\ref{C_ContInV}),
(\ref{C_KirhInV}) at all internal vertices,}
\label{Sp_IntCnd}\\
\left\{\phi,\psi\right\}=0\quad\text{on $\Gamma$}.\label{SpDirBnd}
\end{eqnarray}
\begin{definition}
By $\mathbf{S}$ we denote the spectral problem described by
(\ref{SpEq1ch})--(\ref{SpDirBnd}).
\end{definition}
It is known that the problem $\mathbf{S}$ has a discrete spectrum of
eigenvalues $0<\lambda_1\leqslant \lambda_1\leqslant
\lambda_2\leqslant\ldots, \lambda_k\to\infty$. Corresponding
eigenfunctions $\{\phi,\psi\}$ can be chosen such that they form an
orthonormal basis in $L_2(\Om).$ Indeed, for scalar problems, i.e.
out-of-the-plane displacements in the mechanical context or
conductivity, the spectral behavior has been explored by von
Below~\cite{vBelow1985}, Nicaise~\cite{Nicaise1995} and others. The
in-plane case discussed here has been treated in Lagnese, Leugering
and Schmidt~\cite{LLS1994}. As a matter of fact, the spectrum
consists of strings of eigenvalues coming from individual edges and
strings of eigenvalues associated with the spectrum of the adjacency
matrix.

\subsection{Dynamical settings.}
Along with the spectral, we consider dynamical system, described
by the two-velocity problem on the each edge of the graph:
\begin{eqnarray}
\frac{1}{k_{i1}^2}u_{tt}^i-u_{xx}^i=0,\quad t>0, \ x\in e_i \label{Eq1ch}\\
\frac{1}{k_{i2}^2}w_{tt}^i-w_{xx}^i=0,\quad t>0, \ x\in e_i
\label{Eq2ch}.
\end{eqnarray}
Here the coefficients $k_{i1}$, $k_{i2}$ play the role of
velocities on the edge $e_i$, $i=1,\ldots N$ in the first and
second channels.

By $\mathcal{F}^T_\Gamma=L_2([0,T], \mathbb{R}^{2m})$ we denote
the space of controls acting on the boundary of the tree. For the
element $F\in \mathcal{F}^T_\Gamma$ we write
$$
F=\{f,g\}=\left\{
\begin{pmatrix}
f^i\\
g^i
\end{pmatrix}
\right\}_{i=1}^m,\quad f^i,g^i\in L_2(0,T).
$$
We will deal with the Dirichlet boundary conditions::
\begin{equation}
\label{C_DirBnd} \{ u,w\}=\{f,g\},\quad \text{on $\Gamma\times
[0,T]$}.
\end{equation}
where ${f,g}\in \mathcal{F}^T_\Gamma$.

\begin{definition}
By $\mathbf{D}$ we denote the dynamical problem on the graph
$\Omega$, described by equations on the edges (\ref{Eq1ch}),
(\ref{Eq2ch}) which satisfies compatibility conditions
(\ref{C_ContInV}), (\ref{C_KirhInV}) at all internal vertexes for
any $t>0$, Dirichlet boundary condition (\ref{C_DirBnd}) and zero
initial conditions $\{u(\cdot,0),w(\cdot,0)\}=\{0,0\}$,
$\{u_t(\cdot,0),w_t(\cdot,0)\}=\{0,0\}$.
\end{definition}
It is known  that for any $ T>0, \  \{u,w\} \in
C([0,T];L^2(\Omega))$. See e.g. \cite{LLS1994,LagneseLeugering2004}.

\section{Inverse dynamical and spectral problems. Connection of the inverse data. }

We use the Titchmarsh-Weyl (TW) matrix as the data for the inverse
spectral problem. For the spectral problem on an interval and on
the half line the TW function is a classical object. For the
inverse spectral problem for the wave equation on trees it was
used in \cite{AK,Y}. The general properties of the $M-$operator
for self-adjoint operators are considered in \cite{AmPe'04, BMN,
BMN1}.

Let us choose $\lambda\notin \mathbb{R}$. We define
$\{\phi^1,\psi^1\}$, $\{\phi^2,\psi^2\}$ --- two solutions of
(\ref{SpEq1ch}), (\ref{Sp_IntCnd}) and the following boundary
conditions:
\begin{equation}
\{\phi^1,\psi^1\}=\begin{pmatrix}(0,0)\\
\hdotsfor{1}\\
(1,0)\\
(0,0)
\end{pmatrix},\quad
\{\phi^2,\psi^2\}=\begin{pmatrix}(0,0)\\
\hdotsfor{1}\\
(0,1)\\
(0,0)
\end{pmatrix},\quad \text{on $\Gamma$},
\end{equation}
where nonzero elements are located at the $i-$th row. Then the TW
matrix $\mathbf{M}(\lambda)$ is defined as
$\mathbf{M}(\lambda)=\{M_{ij}(\lambda)\}_{i,j=1}^n$ where each
$M_{ij}(\lambda)$ is a $2\times 2$ matrix defined by
\begin{equation}
M_{ij}(\lambda)=\begin{pmatrix}\phi^1_x(v_j,\lambda) & \psi^1_x(v_j,\lambda)\\
\phi^2_x(v_j,\lambda) & \psi^2_x(v_j,\lambda)
\end{pmatrix},\quad 1\leqslant i,j\leqslant m.
\end{equation}

Let us consider the nonhomogeneous Dirichlet boundary condition
\begin{equation}
\label{SpDirF} \{\phi,\psi\}=\{\zeta,\nu\},
\end{equation}
and let $\{\phi,\psi\}$ be the solution to (\ref{SpEq1ch}),
(\ref{Sp_IntCnd}), (\ref{SpDirF}). The Titchmarsh-Weyl matrix
connects the values of $\{\phi,\psi\}$ on the boundary and the
values of its derivative $\{\phi_x,\psi_x\}$ on the boundary:
\begin{equation}
\label{M_def}
\{\phi_x,\psi_x\}=\mathbf{M}(\lambda)\{\zeta,\nu\},\quad \text{on
$\Gamma$}.
\end{equation}

We set up the \emph{ spectral inverse problem} as follows: given
the TW matrix $\mathbf{M}(\lambda)$, $\lambda\notin \mathbb{R}$,
to recover the graph (lengths of edges, connectivity and angles
between edges) and parameters of the system (\ref{SpEq1ch}), i.e.
the set of coefficients $\{k_{i1},k_{i2}\}_{i=1}^N$.

Let $\{u,w\}$ be the solution to the problem $\mathbf{D}$ with the
boundary control $\{f,g\}\in \mathcal{F}^T_\Gamma$. We introduce
the \emph{dynamical response operator} (the dynamical
Dirichlet-to-Neumann map) to the problem $\mathbf{D}$ by the rule
\begin{equation}
\label{RespOp}
R^T\{f,g\}(t)=\{u_x(\cdot,t),w_x(\cdot,t)\}\Bigl|_\Gamma,\quad
0\leqslant t\leqslant T.
\end{equation}
The response operator has the form of a convolution:
\begin{equation}
\left(R^T\{f,g\}\right)\left(t\right)=\left(\mathbf{R}*\{f,g\}\right)\left(t\right),
\quad 0\leqslant t\leqslant T
\end{equation}
where $\mathbf{R}(t)=\{R_{ij}(t)\}_{i,j=1}^n$ and each
$R_{i,j}(t)$ is a $2\times 2$ matrix. The entries $R_{ij}(t)$ are
defined by the following procedure. We set up two dynamical
problems (\ref{Eq1ch}), (\ref{Eq2ch}), (\ref{C_ContInV}),
(\ref{C_KirhInV}) and the boundary conditions given by
\begin{equation}
\label{Dyn_BC}
\{U^1(\cdot,t),W^1(\cdot,t)\}=\begin{pmatrix}(0,0)\\
\hdotsfor{1}\\
(\delta(t),0)\\
(0,0)
\end{pmatrix},\,\,
\{U^2(\cdot,t),W^2(\cdot,t)\}=\begin{pmatrix}(0,0)\\
\hdotsfor{1}\\
(0,\delta(t))\\
(0,0)
\end{pmatrix},\,\, \text{on $\Gamma$}.
\end{equation}
In the above notations, the only nonzero rows is $i-$th. Then
\begin{equation}
\label{RespDef}
R_{ij}(t)=\begin{pmatrix}U^1_x(v_j,t) & W^1_x(v_j,t)\\
U^2_x(v_j,t) & W^2_x(v_j,t)
\end{pmatrix}.
\end{equation}
So, to construct the entries of $\mathbf{R}$, we need to set up
the boundary condition at $i-$th boundary point in the first and
second channels, while having other boundary points fixed (impose
homogeneous Dirichlet conditions there) and measure the response
at $j-$th boundary point in the first and second channels.

We set up the \emph{ dynamical inverse problem} as follows: given
the response operator $R^T$ (\ref{RespOp}) (or what is equivalent,
the matrix $\mathbf{R}(t)$, $t\in [0,T]$), for large enough $T,$
to recover the graph (lengths of edges, connectivity and angles
between edges) and parameters for the dynamical system
(\ref{Eq1ch}), (\ref{Eq2ch}), i.e. velocities on the edges.

The connection of the spectral and dynamical data is known and was
used for studying the inverse spectral and dynamical problems, see
for example \cite{KKLM,AK,AMR}. Let $\{f,g\}\in
\mathcal{F}^T_\Gamma\cap (C_0^\infty(0,+\infty))^{2m}$ and
$$
\widehat{\{f,g\}}(k):=\int_0^\infty \{f(t),g(t)\}e^{ikt}\,dt
$$
be its Fourier transform. The equations (\ref{Eq1ch}),
(\ref{Eq2ch}) and (\ref{SpEq1ch}) are connected by the Fourier
transformation: going formally in (\ref{Eq1ch}), (\ref{Eq2ch})
over to the Fourier transform, we obtain (\ref{SpEq1ch}) with
$\lambda=k^2$. It is not difficult to check (see, e.g. \cite{AMR,
AK} that the response operator (response function) and
Titchmarsh-Weyl matrix are connected by the same transformation:
\begin{equation}
\label{Sp_Dyn} \mathbf{M}(k^2)=\int_0^\infty
\mathbf{R}(t)e^{ikt}\,dt.
\end{equation}

\section{Solution of the inverse problem. The case of two intervals.}

We start with the inverse problem for two connected intervals. For
the two-velocity system even this simple situation is nontrivial.

Suppose that a tree consists of two edges, $e_1$ and  $e_2$ with
the (unknown) lengths $l_1$ and $l_2$. The angle between edges we
denote by $\alpha:=\alpha_{12}$ and
$\Gamma=\partial\Omega=\{v_1,v_2\}$ and the only internal point is
$v_3$. We consider the dynamical problem $\mathbf{D}$ on the tree
and show that we need to know response operator or TW function
associated with one boundary vertex only to recover the graph.

Let us consider the problem (\ref{Eq1ch}), (\ref{Eq2ch}),
(\ref{C_ContInV}), (\ref{C_KirhInV}) with the boundary conditions
given by
\begin{equation}
\label{BCFirstCh}
\{u,w\}=\begin{pmatrix}(\delta(t), 0)\\
(0, 0)
\end{pmatrix},\quad\text{on $\Gamma$}.
\end{equation}
The solution of the above problem can be evaluated explicitly. For
$0\leqslant t\leqslant \frac{l_1}{k_{11}}$ it is given by
\begin{eqnarray*}
u^1(x,t)=\delta\left(t-\frac{x}{k_{11}}\right),\\
w^1(x,t)=0.
\end{eqnarray*}
On the time interval $\frac{l_1}{k_{11}}<t<
\frac{2l_2}{\min{\{k_{21},k_{22\}}}}$ on the first edge we have
\begin{eqnarray*}
u^1(x,t)=\delta\left(t-\frac{x}{k_{11}}\right)+a_1\delta\left(t+\frac{x}{k_{11}}-\frac{2l_1}{k_{11}}\right),\\
w^1(x,t)=b_1\delta\left(t+\frac{x}{k_{12}}-\gamma_{12}\right),\quad
\gamma_{12}=\frac{l_1}{k_{11}}+\frac{l_1}{k_{12}}.
\end{eqnarray*}
and on the second edge
\begin{eqnarray*}
u^2(x,t)=a_2\delta\left(t-\frac{x}{k_{21}}-\gamma_{21}\right),\quad \gamma_{21}=\frac{l_1}{k_{11}}-\frac{l_1}{k_{21}},\\
w^2(x,t)=b_2\delta\left(t-\frac{x}{k_{22}}-\gamma_{22}\right),\quad
\gamma_{22}=\frac{l_1}{k_{11}}-\frac{l_1}{k_{22}}.
\end{eqnarray*}
In the formulas above, the coefficients $a_{1}$, $a_2$, $b_1$,
$b_2$ are unknown.

>From the condition (\ref{C_ContInV}) we obtain that
\begin{equation}
\label{c1}
\begin{pmatrix}
1+a_1\\
b_1
\end{pmatrix}=
R_\alpha \begin{pmatrix}
a_2\\
b_2
\end{pmatrix}.
\end{equation}
Condition (\ref{C_KirhInV}) implies
\begin{equation}
\label{Kir_2str} D_1'\begin{pmatrix}
-\frac{1}{k_{11}}+\frac{a_1}{k_{11}}\\
\frac{1}{k_{12}}
\end{pmatrix}=R_\alpha D_2'
\begin{pmatrix}
-\frac{a_2}{k_{21}}\\
-\frac{b2}{k_{22}}
\end{pmatrix}.
\end{equation}
After introducing the notation
\begin{equation*}
D_i=\begin{pmatrix}
k_{i1} & 0\\
0 & k_{i2}
\end{pmatrix},\quad i=1,\ldots N,
\end{equation*}
we can rewrite (\ref{Kir_2str}) as
\begin{equation}
\label{c2}
D_1\begin{pmatrix}
-1+a_1\\
b_1
\end{pmatrix}=-R_\alpha D_2
\begin{pmatrix}
a_2\\
b_2
\end{pmatrix}.
\end{equation}
Combining (\ref{c1}) and (\ref{c2}), we obtain
\begin{equation}
\label{c3}
D_1\begin{pmatrix}
-1+a_1\\
b_1
\end{pmatrix}=-R_\alpha D_2 R_{-\alpha}
\begin{pmatrix}
1+a_1\\
b_1
\end{pmatrix}.
\end{equation}
Let us now consider the problem (\ref{Eq1ch}), (\ref{Eq2ch}),
(\ref{C_ContInV}), (\ref{C_KirhInV}) with the following boundary
condition at the second channel:
\begin{equation}
\label{BCSecCh}
\{u,v\}=\begin{pmatrix}(0,\delta(t)) \\
(0, 0)
\end{pmatrix},\quad \text{on $\Gamma$}.
\end{equation}
The computations similar to ones for the case of the boundary
condition at the first channel show that the following condition
must hold:
\begin{equation}
\label{c4} D_1\begin{pmatrix}
\widetilde b_1\\
-1+\widetilde a_1
\end{pmatrix}=-R_\alpha D_2 R_{-\alpha}
\begin{pmatrix}
\widetilde b_1\\
1+\widetilde a_1
\end{pmatrix},
\end{equation}
where the coefficients $\widetilde a_{1}$, $\widetilde a_2$,
$\widetilde b_1$, $\widetilde b_2$ are unknown.

The function $u^1(x,t)$, the first component of the solution to
(\ref{Eq1ch}), (\ref{Eq2ch}), (\ref{C_ContInV}), (\ref{C_KirhInV}),
(\ref{BCFirstCh}) has the following representation on the time
interval
$\frac{2l_1}{k_{11}}<t<\frac{2l_1}{k_{11}}+2\min_{i,j}\{\frac{l_i}{k_{ij}}\}$
\begin{equation}
\label{N}
u^1(x,t)=\delta\left(t-\frac{x}{k_{11}}\right)+a_1\delta\left(t+\frac{x}{k_{11}}-\frac{2l_1}{k_{11}}\right)-
a_1\delta\left(t-\frac{x}{k_{11}}-\frac{2l_1}{k_{11}}\right)
\end{equation}
Thus (see the definition of the response operator (\ref{RespDef}))
for such $t$:
\begin{equation}
\label{r1}
\{R_{11}\}_{11}(t)=u^1_x(0,t)=-\frac{1}{k_{11}}\delta'(t)+\frac{2a_1}{k_{11}}\delta'\left(t-\frac{2l_1}{k_{11}}\right).
\end{equation}
A similar argument shows that for
$\frac{l_1}{k_{11}}+\frac{l_1}{k_{12}}<t<\frac{l_1}{k_{11}}+\frac{l_1}{k_{12}}+
2\min_{i,j}\{\frac{l_i}{k_{ij}}\}$
\begin{equation}
\label{r2} \{R_{11}\}_{12}(t)=w^1_x(0,t)=2
\frac{b_1}{k_{12}}\delta'\left(t-\gamma_{12}\right).
\end{equation}
Thus using the $\{R_{11}\}_{12}$, $\{R_{11}\}_{12}$ components of
the response operator on the described time intervals, we can
determine $k_{11}$, $k_{12}$, $a_1$, $b_1$, $l_1$. Applying the
same argument to the problem (\ref{Eq1ch}), (\ref{Eq2ch}),
(\ref{C_ContInV}), (\ref{C_KirhInV}), (\ref{BCSecCh}), we conclude
that $\{R_{11}\}_{21}$, $\{R_{11}\}_{22}$ components of the
response function determine $\widetilde a_1$, $\widetilde b_1$.

Let us introduce the notations
\begin{eqnarray*}
\xi_1=1+a_1,\quad \eta_1=b_1,\\
\widetilde \xi_1=\widetilde b_1,\quad \widetilde
\eta_1=1+\widetilde a_1,
\end{eqnarray*}
and rewrite (\ref{c3}), (\ref{c4}) as
\begin{eqnarray}
\left(R_\alpha D_2 R_{-\alpha}+D_1\right)\begin{pmatrix}
\xi_1\\
\eta_1
\end{pmatrix}= 2D_1\begin{pmatrix}
1\\
0
\end{pmatrix},\label{c5}\\
\left( R_\alpha D_2 R_{-\alpha}+D_1\right)
\begin{pmatrix}
\widetilde \xi_1\\
\widetilde \eta_1
\end{pmatrix}=2D_1\begin{pmatrix}
0\\
1
\end{pmatrix}\label{c6}.
\end{eqnarray}
In the forward problem, the equations (\ref{c5}) and (\ref{c6}) can
be used for the determination of the reflection and transmission
coefficients $a_1$, $b_1$, $\widetilde a_1$, $\widetilde b_1$. Since
the (given) matrixes $R_\alpha D_2 R_{-\alpha}+D_1$ and $D_1$ are
positive definite (it is easy to check that $R_\alpha D_j
R_{-\alpha}\geqslant \min\{k_{j1},k_{j2}\}I$), these coefficients
are uniquely determined from (\ref{c5}), (\ref{c6}). Moreover, we
necessarily have
\begin{equation}
\begin{pmatrix}
\xi_1\\
\eta_1
\end{pmatrix}\not=\begin{pmatrix}
\widetilde \xi_1\\
\widetilde \eta_1
\end{pmatrix}
\end{equation}
On the other hand, in the inverse problem we know $D_1$, $\xi_1$,
$\eta_1$, $\widetilde \xi_1$, $\widetilde \eta_1$ (we can
determine all these coefficients from $R_{11}$ component of the
response operator). Thus equations (\ref{c5}), (\ref{c6})
determine the matrix $A=R_\alpha D_2 R_{-\alpha}$. Using the
invariants of matrix --- the determinant and trace, we can find
the matrix $D_2$ from the equations
\begin{eqnarray}
k_{21}+k_{22}=\tr{A},\quad k_{21}k_{22}=\det{A}.
\end{eqnarray}
The existence of angle $\alpha$ follows from the spectral theorem
--- $R_\alpha$ ``diagonalize" operator $A$.

Notice that if $D_2=\gm I,$ there is no dependence on $\alpha$ in
equations (\ref{c5}), (\ref{c6}).

We point out that we still need to determine the length of the
second edge. For this aim we could analyze the representation of the
solutions $u^1$, $w^1$ on a sufficiently large  time interval. It
would lead to an increasing number of terms in (\ref{N}) and
(\ref{r1}), (\ref{r2}). Instead of that we will develop the method
which works for general trees following ideas of \cite{AK}. Let us
consider the new tree, consisting of one edge: $\widetilde
\Omega=e_2$. Idea of the method is to recalculate the TW matrix for
the new tree using the TW matrix and response operator for the whole
tree $\Omega$ and the data that we obtained on the first step, i.e.
parameters of the first edge and the angle between edges.

Let $\{U,W\}$ be the solution of (\ref{SpEq1ch}),
(\ref{Sp_IntCnd}) on $\Omega$ with boundary conditions
\begin{equation}
\label{TwoSt1} \{U,W\}=\{\zeta,\nu\},\,\,\text{at $v_1$},\quad
\{U,W\}=\{0,0\}, \,\,\text{at $v_2$}.
\end{equation}
The compatibility conditions at the internal vertex $v_3$ are:
\begin{eqnarray}
\begin{pmatrix}
U^1(v_3,\lambda)\\
W^1(v_3,\lambda)
\end{pmatrix}=R_{\alpha}\begin{pmatrix}
U^2(v_3,\lambda)\\
W^2(v_3,\lambda)
\end{pmatrix},\label{g1}\\
D'_1
\begin{pmatrix}
U^1_x(v_3,\lambda)\\
W^1_x(v_3,\lambda)
\end{pmatrix}=R_{\alpha}D'_2\begin{pmatrix}
U^2_x(v_3,\lambda)\\
W^2_x(v_3,\lambda)
\end{pmatrix}.\label{g2}
\end{eqnarray}
Let $\widetilde M(\lambda)$ be the TW matrix for the tree
$\widetilde\Omega$. We see that $\widetilde M_{11}(\lambda)$ --
component associated with the ``new" boundary point $v_3$
satisfies equation:
\begin{equation}
\begin{pmatrix}
U^2_x(v_3,\lambda)\\
W^2_x(v_3,\lambda)
\end{pmatrix}=\widetilde {M}_{11}(\lambda)\begin{pmatrix}
U^2(v_3,\lambda)\\
W^2(v_3,\lambda)
\end{pmatrix}.\label{TiM}
\end{equation}
>From (\ref{g1})--(\ref{TiM}) it follows
\begin{equation}
\label{m1} D'_1
\begin{pmatrix}
U^1_x(v_3,\lambda)\\
W^1_x(v_3,\lambda)
\end{pmatrix}=R_{\alpha}D'_2\widetilde M_{11}(\lambda)R_{-\alpha}
\begin{pmatrix}
U^1(v_3,\lambda)\\
W^1(v_3,\lambda)
\end{pmatrix}.
\end{equation}
We emphasize that in (\ref{m1}) we know everything but the matrix
$\widetilde M(\lambda)$. Choosing different boundary conditions for
the problem in (\ref{TwoSt1}), we can get linear independent vectors
$\begin{pmatrix}
U^1(v_3,\lambda)\\
W^1(v_3,\lambda)
\end{pmatrix}$
in (\ref{m1}). Thus (\ref{m1}) determines the matrix $\widetilde
M_{11}(\lambda)$. The matrix $\widetilde M_{11}$ uniquely
determines the corresponding component of the response operator
$\widetilde R_{11}$ (see (\ref{Sp_Dyn})). The latter operator in
turn, allows us to find the parameters of the second edge, exactly
as $R_{11}$ determines the parameters of the first edge. (In our
simple case of the two edge tree, the only parameter which we need
to recover is the length of the second edge.)

We conclude the results of the present section in the following
statement
\begin{theorem}
Let $\Omega$ be the tree consisted of two edges. Then the tree and
the parameters of the systems (\ref{Eq1ch}), (\ref{Eq2ch}) and
(\ref{SpEq1ch}) are determined by  the $2 \times 2$ matrix
$M_{11}(\lambda)$ the diagonal element of the TW matrix,
associated with the first boundary point.
\end{theorem}

In the remaining part of the  paper we will deal with the dynamical
approach to the inverse problem. We will prove the following
statement:
\begin{remark}
The  $2 \times 2$ matrix $R_{11}(t)$ for $t \in [0,T]$ where $T>
2d(\Omega)$ determines the tree and the parameters of the systems
(\ref{Eq1ch}), (\ref{Eq2ch}) and (\ref{SpEq1ch})
\end{remark}

\section{Solution of the inverse problem. The case of a star graph.}

Suppose that the tree is a star graph with edges $e_1,\ldots e_n$,
$e_i=[z_i,z_0]$, $i=1,\ldots n$. To recover the graph and the
parameters of the system (\ref{Eq1ch}), (\ref{Eq2ch}) we use the
diagonal elements of the response operator (or diagonal elements
of the Weil matrix), associated to the boundary vertices
$z_1,\ldots z_{n-1}$.

Let us  set up the initial-value problem (\ref{Eq1ch}),
(\ref{Eq2ch}), (\ref{C_ContInV}), (\ref{C_KirhInV}) with the
boundary conditions given by the first equation in (\ref{Dyn_BC}).
We suppose that the only nonzero boundary conditions are given at
the first channel of the $i-th$ boundary vertex, $i\not=n$.
Analyzing the solution $\{u,w\}$ of these problems in a way we did
for the case of the graph of two edges, we obtain that on the time
interval $\frac{l_i}{k_{i1}}<t< \frac{l_i}{k_{i1}} +2 \min_{i,j}
\{\frac{l_i}{k_{ij}}\}$ on the $i-$th edge we have
\begin{eqnarray*}
u^i(x,t)=\delta\left(t-\frac{x}{k_{i1}}\right)+a^i_i\delta\left(t+\frac{x}{k_{i1}}-\frac{2l_i}{k_{i1}}\right),\\
w^i(x,t)=b^i_i\delta\left(t+\frac{x}{k_{i2}}-\gamma_{i2}\right),\quad
\gamma_{i2}=\frac{l_i}{k_{i1}}+\frac{l_i}{k_{i2}}.
\end{eqnarray*}
and on other edges ($j=1\ldots n$, $j\not=i$):
\begin{eqnarray*}
u^j(x,t)=a^i_j\delta\left(t-\frac{x}{k_{j1}}-\gamma_{j1}\right),\quad \gamma_{j1}=\frac{l_i}{k_{i1}}-\frac{l_i}{k_{j1}},\\
w^j(x,t)=b^i_j\delta\left(t-\frac{x}{k_{j2}}-\gamma_{j2}\right),\quad
\gamma_{j2}=\frac{l_i}{k_{i1}}-\frac{l_i}{k_{j2}}.
\end{eqnarray*}
Where $a^i_j$, $b^i_j$, $i=1,\ldots,n-1$, $j=1,\ldots,n$ are
reflection and transmission coefficients associated with the
$i-th$ vertex. Let us introduce new parameters
\begin{eqnarray}
\xi_i=1+a^i_i,\quad \eta_i=b^i_i, \quad i=1,\ldots, n-1.
\end{eqnarray}
The compatibility conditions (\ref{C_ContInV}), (\ref{C_KirhInV})
at the internal vertex $z_0$ (we need to rewrite them in a way we
did for the case of two edges) lead to the following equalities
(cf. (\ref{c5})):
\begin{eqnarray}
\left(\sum_{j\not= i}R_{\alpha_{ij}} D_j
R_{-\alpha_{ij}}+D_i\right)\begin{pmatrix}
\xi_i\\
\eta_i
\end{pmatrix}=2D_i
\begin{pmatrix}
1\\
0
\end{pmatrix},\quad i=1,\ldots n-1\label{d1}.
\end{eqnarray}
Let us now set up the initial-value problem with the delta
function in the second channel at $i-$th boundary point,
$i\not=n$, which is given by (\ref{Eq1ch}), (\ref{Eq2ch}),
(\ref{C_ContInV}), (\ref{C_KirhInV}) and the boundary conditions
given by the second equation in (\ref{Dyn_BC}). We can obtain and
analyze the representation for the solutions $\{\widetilde
u,\widetilde w\}$ of these problems. Let $\widetilde a_i$,
$\widetilde b_i$, $i=1,\ldots,n-1$, $j=1,\ldots,n$ be the
reflection and transmission coefficients. Introducing  new
parameters
\begin{eqnarray}
\widetilde \xi_i=\widetilde b^i_i,\quad \widetilde
\eta_i=1+\widetilde a^i_i,\quad i=1,\ldots n-1
\end{eqnarray}
and making use of the compatibility conditions (\ref{C_ContInV}),
(\ref{C_KirhInV}) at the internal vertex $z_0$, we obtain the
following equalities (cf. (\ref{c6})):
\begin{eqnarray}
\left(\sum_{j\not= i}R_{\alpha_{ij}} D_j
R_{-\alpha_{ij}}+D_i\right)\begin{pmatrix}
\widetilde \xi_i\\
\widetilde \eta_i
\end{pmatrix}=2D_i
\begin{pmatrix}
0\\
1
\end{pmatrix},\quad i=1,\ldots n-1.\label{d2}
\end{eqnarray}
The matrices $\left(\sum_{j\not= i}R_{\alpha_{ij}} D_j
R_{-\alpha_{ij}}+D_i\right)$, $i=1,\ldots ,n-1$ are positive
definite. If all angles between edges and all matrices $D_j$ are
known, the systems (\ref{d1}), (\ref{d2}) can be solved for $\xi_i$,
$\eta_i$, $\widetilde \xi_i$, $\widetilde \eta_i$. Note that
necessarily
$$
\begin{pmatrix}
\widetilde \xi_i\\
\widetilde \eta_i
\end{pmatrix}\not=\begin{pmatrix}
\xi_i\\
\eta_i
\end{pmatrix}, \quad i=1,\ldots n-1.
$$
In the situation of the inverse problem, using the diagonal elements
$\{R_{ii}\}$, $i=1,\ldots,n-1$ of the response operator, we can
determine the reflection and transmission coefficients $a^i_i$,
$b^i_i$, $\widetilde a^i_i$, $\widetilde b^i_i$, as well as $l_i$,
$D_i$ for $i=1,\ldots ,n-1$. Indeed: analyzing the solution to the
dynamical system $\mathbf{D}$ with the boundary condition given by
the delta function in the first channel at the $i-$th boundary
vertex, it is easy to see (cf. (\ref{r1}), (\ref{r2})) that:
\begin{eqnarray*}
\{R_{ii}\}_{11}(t)={u^i}_x(0,t)=-\frac{1}{k_{i1}}\delta'(t)+\frac{2a^i_i}{k_{i1}}\delta'\left(t-\frac{2l_i}{k_{i1}}\right),
\\ \quad
\frac{2l_i}{k_{i1}}<t<\frac{2l_i}{k_{i1}}+2\min_{i,j}\{\frac{l_i}{k_{ij}}\},\\
\{R_{ii}\}_{12}(t)={w^i}_x(0,t)=\frac{b^i_i}{k_{i2}}\delta'\left(t-\gamma_{i2}\right),
\\ \quad
\gamma_{i2}<t<\gamma_{i2}+2\min_{i,j}\{\frac{l_i}{k_{ij}}\}.
\end{eqnarray*}
The above representation allows one to determine $a^i_i$, $b^i_i$,
$l_i$, $k_{i1}$ for $i=1,\ldots,n-1$. Analyzing the solutions to
the dynamical system $\mathbf{D}$ with the boundary condition
given by the delta function at the second channel of the $i-$th
boundary vertex, we can determine $\widetilde a^i_i$, $\widetilde
b^i_i$, $k_{i2}$ for $i=1,\ldots,n-1$.

Thus, since the vectors $\begin{pmatrix}
\widetilde \xi_i\\
\widetilde \eta_i
\end{pmatrix}$ and $\begin{pmatrix}
\xi_i\\
\eta_i
\end{pmatrix}$
in (\ref{d1}), (\ref{d2}) are known and necessarily different,
equations (\ref{d1}), (\ref{d2}) completely determines the
matrixes $A_i$,
\begin{equation}
\label{d3} A_i=\sum_{j\not= i}R_{\alpha_{ij}} D_j
R_{-\alpha_{ij}}+D_i, \quad i=1\ldots n.
\end{equation}
In (\ref{d3}) we do not know $n-1$ angles between edges and the
matrix $D_n$. We  use the system (\ref{d3}) to determine them in
the same way we did it for the case of two intervals, but
calculations are more involved.

Let us consider the condition (\ref{d3}) for $i=k$ and for $i=l$:
\begin{eqnarray}
D_k+R_{\alpha_{kl}} D_l R_{-\alpha_{kl}}+\sum_{j\not= k,\,l}R_{\alpha_{kj}} D_j R_{-\alpha_{kj}}=A_k, \label{S1}\\
D_{l}+R_{\alpha_{lk}} D_k R_{-\alpha_{lk}}+\sum_{j\not=
k,\,l}R_{\alpha_{lj}} D_j R_{-\alpha_{lj}}=A_l,\label{S2}
\end{eqnarray}
where the matrices $A_k$ and $A_l$ are known. Note that after the
multiplication of (\ref{S2}) by $R_{\alpha_{kl}}$ from the left
and by $R_{-\alpha_{kl}}$ from the right, and using that
$R_{\alpha_{kl}}=R_{-\alpha_{lk}}$,
$R_{\alpha_{kl}}R_{\alpha_{lj}}=R_{\alpha_{kj}}$, we obtain
\begin{equation*}
A_k=R_{\alpha_{kl}}A_lR_{-\alpha_{kl}}.
\end{equation*}
The angle $\alpha_{kl}$ can now be found using the spectral
theorem. Repeating this procedure for various $i$, $l$ we can
determine all angles. After that we can use any of the conditions
(\ref{d3}) to determine $D_n$. Indeed, taking $i=k$ we have
\begin{equation}
R_{\alpha_{kn}} D_n R_{-\alpha_{kn}}=B_k
\end{equation}
for some known matrix $B_k$. Then
\begin{eqnarray}
k_{n1}+k_{n2}=\tr{B_k},\quad k_{n1}k_{n2}=\det{B_k}.
\end{eqnarray}

The next step is crucial for solving the inverse problem: we have
already recovered a part of the tree, and our next aim is to find
the inverse data for the smaller ``new" tree, using the initial
inverse data and information that we obtained on the previous steps.

Let us consider the new tree, consisting of the one edge
$\widetilde\Omega=e_n=[z_0,z_n]$. By $\{\Phi,\Psi\}$ we denote the
solution to (\ref{SpEq1ch}), (\ref{Sp_IntCnd}) and the following
boundary conditions
\begin{equation}
\label{TwoSt10} \{\Phi,\Psi\}=\{\zeta,\nu\},\,\,\text{at
$z_1$},\quad \{\Phi,\Psi\}=\{0,0\}, \,\,\text{at $z_i,\,\,
2\leqslant i \leqslant n$}.
\end{equation}
As in the case of  a graph consisting of two intervals, our goal is
to obtain the coefficient $\widetilde M_{11}$ of the TW-matrix for
$\widetilde\Omega$, associated with the ``new" boundary edge $z_0$.
Note that we can assume that we have already recovered the
information about all other edges and angles between them. So we
have in hands the matrices $D_i'$, and $\alpha_{i,n}$ for
$i=1,\ldots n$.

Note that solution to \ref{SpEq1ch}), (\ref{Sp_IntCnd}),
(\ref{TwoSt10}) on the edge $e_1$ solves the Cauchy problem
\begin{eqnarray}
-k_{11}^2\Phi_{xx}^1=\lambda \Phi^1,\quad
-k_{12}^2\Psi_{xx}^1=\lambda \Psi^1,\quad x\in e_1 \\
\{\Phi^1(z_1),\Psi^1(z_1)\}= \{\zeta,\nu\},\\
\begin{pmatrix}
\Phi^1_x(z_1) \\
\Psi^1_x(z_1)
\end{pmatrix}=
\{M_{11}(\lambda)\}
\begin{pmatrix}
\zeta \\
\nu
\end{pmatrix},
\end{eqnarray}
and on edges $e_2,\ldots,e_{n-1}$ solves the Cauchy problems
\begin{eqnarray}
-k_{i1}^2\Phi_{xx}^i=\lambda \Phi^i,\quad
-k_{i2}^2\Psi_{xx}^i=\lambda \Psi^i,\quad x\in e_i \\
\{\Phi^1(z_i),\Psi^1(z_i)\}= \{0,0\},\\
\begin{pmatrix}
\Phi^i_x(z_i) \\
\Psi^i_x(z_i)
\end{pmatrix}=
\{M_{1i}(\lambda)\}
\begin{pmatrix}
\zeta \\
\nu
\end{pmatrix}.
\end{eqnarray}
Thus the function $\{\Phi,\Psi\}$ and its derivative is known on the
edges $e_1,\ldots,e_{n-1}$. At the internal vertex $z_0$
compatibility conditions hold:
\begin{eqnarray}
\begin{pmatrix}
\Phi^1(z_0,\lambda)\\
\Psi^1(z_0,\lambda)
\end{pmatrix}=R_{\alpha_{in}}\begin{pmatrix}
\Phi^n(z_0,\lambda)\\
\Psi^n(z_0,\lambda)
\end{pmatrix},\\
D'_1
\begin{pmatrix}
\Phi^1_x(z_0,\lambda)\\
\Psi^1_x(z_0,\lambda)
\end{pmatrix}=\sum_{j\not=1,n}R_{\alpha_{1j}}D'_j\begin{pmatrix}
\Phi^j_x(z_0,\lambda)\\
\Psi^j_x(z_0,\lambda)
\end{pmatrix}+R_{\alpha_{1n}}D'_n\begin{pmatrix}
\Phi^n_x(z_0,\lambda)\\
\Psi^n_x(z_0,\lambda)
\end{pmatrix}.
\end{eqnarray}
Using these conditions and the definition of the component of
TW-matrix associated with the $n-$th edge:
\begin{equation}
\begin{pmatrix}
\Phi^n_x(z_0,\lambda)\\
\Psi^n_x(z_0,\lambda)
\end{pmatrix}=\widetilde M_{11}(\lambda)\begin{pmatrix}
\Phi^n(z_0,\lambda)\\
\Psi^n(z_0,\lambda)
\end{pmatrix},
\end{equation}
we get the equations
\begin{eqnarray}
\label{m2} D'_i
\begin{pmatrix}
\Phi^1_x(z_0,\lambda)\\
\Psi^1_x(z_0,\lambda)
\end{pmatrix}=\sum_{j\not=1,n} R_{\alpha_{1j}}D'_j \begin{pmatrix}
\Phi^j_x(z_0,\lambda)\\
\Psi^j_x(z_0,\lambda)
\end{pmatrix}+ \\
R_{\alpha_{1n}}D'_n\widetilde M_{11}(\lambda)R_{-\alpha_{1n}}
\begin{pmatrix}
\Phi^1(z_0,\lambda)\\
\Psi^1(z_0,\lambda)
\end{pmatrix}\notag.
\end{eqnarray}
Choosing the different boundary conditions at the $i-$th boundary
point, we can get vectors $\begin{pmatrix}
\Phi^i(z_0,\lambda)\\
\Psi^i(z_0,\lambda)
\end{pmatrix}$ in (\ref{m2}) to be linearly independent. Since we know all other
terms in (\ref{m2}), this equation determines $\widetilde
M_{11}(\lambda)$. Using the connection of the dynamical and
spectral data (\ref{Sp_Dyn}), we can recover the $\widetilde
R_{11}$ component of the response function associated with
$\widetilde\Omega$ and reduce our problem to the inverse problem
for one edge. (Really we still need to recover only the length of
the $n-$th edge.)

We combine all results of this section in
\begin{theorem}
Let $\Omega$ be the a star graph consisted of $n$ edges. Then the
graph and the parameters of the systems (\ref{Eq1ch}),
(\ref{Eq2ch}) and (\ref{SpEq1ch}), are determined by the diagonal
elements ($2\times 2$ matrices) $M_{ii}(\lambda)$ $1\leqslant
i\leqslant n-1$ of the $TM$ matrix.
\end{theorem}

We proceed to prove the following statement:
\begin{remark}
The diagonal elements ($2\times 2$ matrices) $R_{ii}(t)$
$1\leqslant i\leqslant n-1$ of the response function for $t
>d(\Omega)$ determine the tree and the parameters of the systems
(\ref{Eq1ch}), (\ref{Eq2ch}) and (\ref{SpEq1ch})
\end{remark}

\section{Solution of the inverse problem. The case of an
arbitrary tree.}

Let $\Omega$ be a finite tree with $m$ boundary   points
$\Gamma=\{\gamma_1,\ldots,\gamma_m\} $. Without  loss of generality
we can assume that the boundary vertex $\gamma_m$ is a root of the
tree. We consider the dynamical problem $\mathbf{D}$ and the
spectral problem $\mathbf{S}$ on $\Omega$. Then the response
function $\mathbf{R}(t)=\{R_{ij}(t)\}_{i,j=1}^{m-1}$ and the TW
matrix $\mathbf{M}(\lambda)=\{M_{ij}(\lambda)\}_{i,j=1}^{m-1}$
associated with all other boundary points are constructed in the
same way as in the section 3.

Let us take two boundary edges, $e_i$ with the length $l_i$ and
velocities in channels $k_{i1}$, $k_{i2}$ and $e_j$ with the
length $l_j$ and velocities in channels $k_{j1}$, $k_{j2}$. These
two edges have one common point if and only if
\begin{equation}
\{R_{ij}\}_{11}(t)=\left\{\begin{array}l =0,\quad \text{for $t<\frac{l_{i1}}{k_{i1}}+\frac{l_{j1}}{k_{j1}}$}\\
\not=0,\quad \text{for
$t>\frac{l_{i1}}{k_{i1}}+\frac{l_{j1}}{k_{j1}}$}
\end{array}\right.,\quad 1\leqslant i,j\leqslant m-1.
\end{equation}
Note that one can use other components of $R_{ij}$ to determine
the connectivity of edges.

This method allows us to divide the boundary edges into groups, such
that edges from one group have a common vertex. Let us take the
first of such  groups, say
$e_1=[z_1,z_0],\ldots,e_{m_0}=[z_{m_0},z_0]$ with boundary vertices
$\gamma_1,\ldots,\gamma_{m_0}$. These edges together with another
edge $e=[z_0,z_{m_0'}]$ form a star graph, the subgraph of $\Omega$.
Note, that using the diagonal elements of the response operator
$\{R_{ii}(t)\}$ or the diagonal elements of the TW-matrix
$M(\lambda)$, $i=1,\ldots,m_0$ by the same method as in the case of
star graph, we can determine angles and velocities for all edges
$e_1,\ldots,e_{m_0}, [z_0,z_{m_0'}]$ and lengths all boundary edges
$e_1,\ldots,e_{m_0}$.

We take the new tree $\widetilde\Omega=\Omega\backslash
\bigcup_{i=1}^{m_0} e_i$. Our goal as in the previous cases is to
calculate $\widetilde M(\lambda)$, the TW-matrix associated with
$\widetilde\Omega$.

By $\{\Phi,\Psi\}$ we denote the solution to (\ref{SpEq1ch}),
(\ref{Sp_IntCnd}) and the following boundary conditions
\begin{equation}
\label{TwoStGen} \{\Phi,\Psi\}=\{\zeta,\nu\},\,\,\text{at
$\gamma_1$},\quad \{\Phi,\Psi\}=\{0,0\}, \,\,\text{at
$\gamma_i,\,\, 2\leqslant i \leqslant m$}.
\end{equation}
Note that the solution to (\ref{SpEq1ch}), (\ref{Sp_IntCnd}),
(\ref{TwoStGen}) on the edge $e_1$ solves the Cauchy problem
\begin{eqnarray}
-k_{11}^2\Phi_{xx}^1=\lambda \Phi^1,\quad
-k_{12}^2\Psi_{xx}^1=\lambda \Psi^1,\quad x\in e_1 \\
\{\Phi^1(z_1),\Psi^1(z_1)\}= \{\zeta,\nu\},\\
\begin{pmatrix}
\Phi^1_x(z_1) \\
\Psi^1_x(z_1)
\end{pmatrix}=
\{M_{11}(\lambda)\}
\begin{pmatrix}
\zeta \\
\nu
\end{pmatrix},
\end{eqnarray}
and on the edges $e_2,\ldots,e_{m_0}$ solves
\begin{eqnarray}
-k_{i1}^2\Phi_{xx}^i=\lambda \Phi^i,\quad
-k_{i2}^2\Psi_{xx}^i=\lambda \Psi^i,\quad x\in e_i \\
\{\Phi^i(z_i),\Psi^i(z_i)\}= \{0,0\},\\
\begin{pmatrix}
\Phi^i_x(z_i) \\
\Psi^i_x(z_i)
\end{pmatrix}=
\{M_{1i}(\lambda)\}
\begin{pmatrix}
\zeta \\
\nu
\end{pmatrix}.
\end{eqnarray}
Thus the function $\{\Phi,\Psi\}$ and its derivative is known on
edges $e_1,\ldots,e_{m_0}$.

At the internal vertex $z_0$ compatibility conditions hold:
\begin{eqnarray}
\begin{pmatrix}
\Phi^1(z_0,\lambda)\\
\Psi^1(z_0,\lambda)
\end{pmatrix}=R_{\alpha_{1m_0'}}\begin{pmatrix}
\Phi^{m_0'}(z_0,\lambda)\\
\Psi^{m_0'}(z_0,\lambda)
\end{pmatrix},\label{Con1}\\
D'_1
\begin{pmatrix}
\Phi^1_x(z_0,\lambda)\\
\Psi^1_x(z_0,\lambda)
\end{pmatrix}=\sum_{j\not=1,m_0'}R_{\alpha_{1j}}D'_j\begin{pmatrix}
\Phi^j_x(z_0,\lambda)\\
\Psi^j_x(z_0,\lambda)
\end{pmatrix}+\label{Con2}\\
R_{\alpha_{1m_0'}}D'_{m_0'}\begin{pmatrix}
\Phi^{m_0'}_x(z_0,\lambda)\\
\Psi^{m_0'}_x(z_0,\lambda)
\end{pmatrix}.\notag
\end{eqnarray}
Using these conditions and the definition of the TW-matrix
associated with the $m_0'-$th edge:
\begin{equation}
\begin{pmatrix}
\Phi^{m_0'}_x(z_0,\lambda)\\
\Psi^{m_0'}_x(z_0,\lambda)
\end{pmatrix}=\widetilde M_{{m_0'}{m_0'}}(\lambda)\begin{pmatrix}
\Phi^{m_0'}(z_0,\lambda)\\
\Psi^{m_0'}(z_0,\lambda)
\end{pmatrix},
\end{equation}
we obtain that
\begin{eqnarray}
\label{m3} D'_1
\begin{pmatrix}
\Phi^1_x(z_0,\lambda)\\
\Psi^1_x(z_0,\lambda)
\end{pmatrix}=\sum_{j=2}^{m_0'} R_{\alpha_{1j}}D'_j \begin{pmatrix}
\Phi^j_x(z_0,\lambda)\\
\Psi^j_x(z_0,\lambda)
\end{pmatrix}+ \\
R_{\alpha_{1{m_0'}}}D'_{m_0'}\widetilde
M_{{m_0'}{m_0'}}(\lambda)R_{-\alpha_{1{m_0'}}}
\begin{pmatrix}
\Phi^1(z_0,\lambda)\\
\Psi^1(z_0,\lambda)
\end{pmatrix}\notag.
\end{eqnarray}
Equation (\ref{m3}) determines the matrix $\widetilde
M_{{m_0'}{m_0'}}(\lambda)$. By definition of the TW-matrix we have
\begin{equation}
\{M_{1j}(\lambda)\}\begin{pmatrix} \zeta\\
\nu
\end{pmatrix}=
\begin{pmatrix}
\Phi^j(e_j)\\
\Psi^j(e_j)
\end{pmatrix},\quad m_0<j<m.
\end{equation}
On the other hand, by the definition for the TW-matrix for the new
tree $\widetilde\Omega$,
\begin{equation}
\begin{pmatrix}
\Phi^j(e_j)\\
\Psi^j(e_j)
\end{pmatrix}=
\{\widetilde M_{m_0'j}(\lambda)\}\begin{pmatrix} \Phi^{m_0'}(z_0)\\
\Psi^{m_0'}(z_0)
\end{pmatrix},\quad m_0<j<m.
\end{equation}
Thus $\{\widetilde M_{m_0'j}(\lambda)\}$ component of the TW
matrix can be found from the equation
\begin{equation}
\{\widetilde M_{m_0'j}(\lambda)\}
\begin{pmatrix}
\Phi^{m_0'}(z_0)\\
\Psi^{m_0'}(z_0)
\end{pmatrix}=
\{M_{1j}(\lambda)\}\begin{pmatrix} \zeta\\
\nu
\end{pmatrix},\quad m_0<j<m.
\end{equation}

To find the components $\widetilde M_{im_0'}(\lambda)$, $m_0<i<m$,
we fix $\gamma_i$ and denote by $\{\Phi,\Psi\}$ the solution to
(\ref{SpEq1ch}), (\ref{Sp_IntCnd}) with the boundary conditions
\begin{equation}
\label{TwoStGen1} \{\Phi,\Psi\}=\{\zeta,\nu\},\,\,\text{at
$\gamma_i$},\quad \{\Phi,\Psi\}=\{0,0\}, \,\,\text{at
$\gamma_j,\,\, j\not= i $}.
\end{equation}
Note that on the edges $e_1,\ldots,e_{m_0}$ $\{\Phi,\Psi\}$
satisfies the equations
\begin{eqnarray}
-k_{j1}^2\Phi_{xx}^j=\lambda \Phi^j,\quad
-k_{j2}^2\Psi_{xx}^j=\lambda \Psi^j,\quad x\in e_j \\
\{\Phi^j(z_j),\Psi^j(z_j)\}= \{0,0\},\\
\begin{pmatrix}
\Phi^j_x(z_j) \\
\Psi^j_x(z_j)
\end{pmatrix}=
\{M_{ij}(\lambda)\}
\begin{pmatrix}
\zeta \\
\nu
\end{pmatrix}.
\end{eqnarray}
Thus, the function $\{\Phi,\Psi\}$ and its derivative are known on
the edges $e_1,\ldots,e_{m_0}$. Using the compatibility conditions
at the internal vertex $z_0$, for every $\begin{pmatrix}\zeta
\\ \nu\end{pmatrix}$ we can find the vectors $\begin{pmatrix}
\Phi^{m_0'}(z_0,\lambda)\\
\Psi^{m_0'}(z_0,\lambda)
\end{pmatrix}
$, $\begin{pmatrix}
\Phi^{m_0'}_x(z_0,\lambda)\\
\Psi^{m_0'}_x(z_0,\lambda)
\end{pmatrix}$.
We emphasize that the function $\{\Phi,\Psi\}$ does not satisfy
zero Dirichlet conditions at $z_0$. Components $\widetilde
M_{im_0'}(\lambda)$, $i=m_0+1,\ldots,m$ can be obtained from the
equations
\begin{equation}
\begin{pmatrix}
\Phi^{m_0'}_x(z_0,\lambda)\\
\Psi^{m_0'}_x(z_0,\lambda)
\end{pmatrix}-\widetilde M_{m_0'm_0'}(\lambda)\begin{pmatrix}
\Phi^{m_0'}(z_0,\lambda)\\
\Psi^{m_0'}(z_0,\lambda)
\end{pmatrix}=\widetilde M_{i{m_0'}}(\lambda)\begin{pmatrix}
\zeta\\
\nu
\end{pmatrix},
\end{equation}

The procedure described reduces the initial problem to the inverse
problem on the smaller subgraph. By repeating the steps a sufficient
amount of times we recover the whole graph and all parameters. We
conclude this section with
\begin{theorem}
Let $\Omega$ be the an arbitrary tree. Then the tree and the
parameters of the systems (\ref{Eq1ch}), (\ref{Eq2ch}) and
(\ref{SpEq1ch}), are determined by the elements ($2\times 2$
matrices) $M_{ij}(\lambda)$ $1\leqslant i,j\leqslant n-1$ of the
TW matrix.
\end{theorem}

 We also infer the following
statement:
\begin{remark}
Let $\Omega$ be the an arbitrary tree. Then the tree and the
parameters of the systems (\ref{Eq1ch}), (\ref{Eq2ch}) and
(\ref{SpEq1ch}), are determined by the elements ($2\times 2$
matrices) $R_{ij}(t)$ $1\leqslant i,j\leqslant m-1$ of the
response matrix for $t>d(\Omega)$.
\end{remark}

\section{Conclusion and further work and open problems}
We have investigated the inverse problem of recovering the material
properties and the geometrical determinants of the topology of a
tree constituted by linear elastic two-velocity channels or
in-plane-models of elastic strings. For a rooted tree this problem
can be solved using measurements at all leaves (besides the root).
The most remarkable novelty is the detection of the angles between
two consecutive elements.

Problems of the same type with varying coefficients will be treated
in a forthcoming publication. Also the numerical recovery of these
quantities is on its way. Clearly, problems involving frames of
Euler-Benroulli and Timoshenko-beams are of great importance too.
Moreover, the question as to which of the properties are detectable,
if the network contains a circuit appears to be an open problem.

\section{Acknowledgments}
We acknowledge the support of the Elite-Network of Bavaria (ENB) and
the DFG-Cluster of Excellence: Engineering of Advanced Materials.

\end{document}